\documentclass[12pt]{article}%
\usepackage{amsfonts}
\usepackage{graphicx}
\usepackage{amsmath}
\usepackage{amssymb}%
\providecommand{\U}[1]{\protect\rule{.1in}{.1in}}
\newtheorem{theorem}{Theorem}

\newtheorem{lemma}[theorem]{Lemma}

\begin{document}

\begin{center}

{\Large On the CLT for stationary Markov chains with trivial tail sigma field}

\bigskip

Magda Peligrad

\bigskip

Department of Mathematical Sciences, University of Cincinnati, PO Box 210025,
Cincinnati, Oh 45221-0025, USA. \texttt{ }
\end{center}

email: peligrm@ucmail.uc.edu\bigskip

\begin{center}
Abstract

\bigskip
\end{center}

In this paper we consider stationary Markov chains with trivial two-sided tail
sigma field, and prove that additive functionals satisfy the central limit
theorem provided the variance of partial sums divided by $n$ is bounded.

\noindent\textit{Keywords:} Markov chains, central limit theorem, tail sigma
field, additive functionals.

\smallskip

\noindent\textit{Mathematical Subject Classification} (2010): 60F05, 60J05, 60G10.

\bigskip

\section{Introduction}

\qquad One of the most useful theorems for stationary sequences is the central
limit theorem for partial sums $S_{n}$ with the normalization $\sqrt{n}$. For
several classes of additive functionals of stationary Markov chains the size
of the variance of partial sums determine the limiting distribution. For
instance for additive functionals of reversible, stationary and ergodic Markov
chains, with centered and square integrable variables, Kipnis and Varadhan
(1986) proved that if $E(S_{n}^{2})/n$ converges to a finite limit, then the
CLT holds. On the other hand, for additive functionals of Harris recurrent and
aperiodic Markov chains with centered and square integrable variables, Chen
(1999, Theorem II. 3.1) proved that if $S_{n}/\sqrt{n}$ is stochastically
bounded, it satisfies the CLT.

These results suggest and motivate the study of limiting distribution for
stationary Markov chains with additive functionals satisfying $\sup_{n}%
E(S_{n}^{2})/n<\infty$. Recently, Peligrad (2020) introduced a new idea, which
involves conditioning with respect to both the past and the future of the
process. By using this approach she proved that functions of a Markov chain
which is stationary and totally ergodic (in the ergodic theoretical sense),
and with $\sup_{n}E(S_{n}^{2})/n<\infty,$ satisfy the CLT,\ provided that a
random centering is used. In this paper we show that the random centering is
not needed if the two-sided sigma field of the Markov chain is trivial.

Our result simply states that if a stationary Markov chain has two-sided tail
sigma field trivial, then any additive functional with finite second moment,
centered at expectation and with $\sup_{n}E(S_{n}^{2})/n<\infty$ satisfies the
central limit theorem (CLT). The interest of such a result consists in the
fact that does not require fine computations of the rate of convergence of
mixing coefficients.

Examples of stationary processes with trivial two-sided tail sigma field
include absolutely regular Markov chains and interlaced mixing Markov chains.
The definitions will be given in this paper. We also refer to Subsection 2.5
in Bradley (2005) for a survey and Bradley (2007) for the proofs of the
results in that survey.

It should be noted that, for a stationary Markov chain, the condition
$\sup_{n}E(S_{n}^{2})/n<\infty$ alone is not enough for CLT (see for instance
Bradley (1989) or Cuny and Lin (2016), Prop. 9.5(ii), among other examples).
On the other hand, if the stationary sequence is not Markov, the conditions
$\sup_{n}E(S_{n}^{2})/n<\infty$ together with the two-sided tail sigma field
is trivial, are not enough for the CLT. Indeed, Bradley (2010) constructed a
stationary sequence, such that any $5$ variables are independent, $\sup
_{n}E(S_{n}^{2})/n<\infty,$ the two-tail sigma field is trivial, but the CLT
does not hold.

\section{Results}

\qquad We assume that $(\xi_{n})_{n\in\mathbb{Z}}$ is a stationary Markov
chain, defined on a complete probability space $(\Omega,\mathcal{F},P)$ with
values in a Polish space $(S,\mathcal{A})$. Denote by $\mathcal{F}_{n}%
=\sigma(\xi_{k},k\leq n)$ and by $\mathcal{F}^{n}=\sigma(\xi_{k},k\geq
n)\,$\ completed with the sets of measure $0$ with respect to $P$. The
marginal distribution on $\mathcal{A}$ is denoted by $\pi(A)=\mathbb{P}%
(\xi_{0}\in A)$.{ }We shall construct the Markov chain in a cannonical way on
$\Omega=S^{Z}$ from a kernel $P(x,A)$, and we assume that an invariant
distribution $\pi$ exists.{ }

We define the two-sided tail sigma field by
\[
\mathcal{T}_{d}=\cap_{n\geq1}(\mathcal{F}_{-n}\vee\mathcal{F}^{n}).
\]
We say that $\mathcal{T}_{d}$ is trivial if for any $A\in\mathcal{T}_{d}$ we
have $P(A)=0$ or $1$. Note that this sigma field might by larger than the
sigma algebra generated by the union of one-sided tail sigma fields defined as
$\mathcal{T}_{l}=\cap_{n\geq1}(\mathcal{F}_{-n})$ and $\mathcal{T}_{r}%
=\cap_{n\geq1}(\mathcal{F}^{n})$. For simplicity, when we refer to the tail
sigma field we shall always understand the two-sided one, $\mathcal{T}_{d}.$

{Let $\mathbb{L}_{0}^{2}(\pi)$ be the set of measurable functions on $S$ such
that $\int f^{2}d\pi<\infty$ and $\int fd\pi=0.$ For a function} ${f}\in
${$\mathbb{L}_{0}^{2}(\pi)$ let }%
\begin{equation}
{X_{i}=f(\xi_{i}),\ S_{n}=\sum\nolimits_{i=1}^{n}X_{i}.} \label{def X,S}%
\end{equation}
{ We denote by }${{||X||}}$ the norm in {$\mathbb{L}^{2}$}$(\Omega
,\mathcal{F},\mathbb{P})$ and $\Rightarrow$ denotes convergence in
distribution.{ }

The goal of this note is to establish the following two results.

\begin{theorem}
\label{ThCLT}Let $(X_{n})_{n\in Z}$ and $(S_{n})_{n\geq1}$ as defined in
(\ref{def X,S}) with
\begin{equation}
\sup_{n}\frac{E(S_{n}^{2})}{n}<\infty. \label{supSn}%
\end{equation}
Assume $(\xi_{n})$\ has trivial tail sigma field $\mathcal{T}_{d}$. Then, for
some $c>0,$ the following limit exists
\[
\lim_{n\rightarrow\infty}\frac{E(|S_{n}|)}{\sqrt{n}}=\frac{c}{\sqrt{2\pi}}%
\geq0\text{ },
\]
and%
\[
\frac{S_{n}}{\sqrt{n}}\Rightarrow N(0,c^{2})\text{ as }n\rightarrow\infty.
\]

\end{theorem}

As a consequence of this result we obtain a necessary and sufficient condition
for the CLT\ for additive functionals of Markov chains with trivial tail sigma field.

\begin{theorem}
\label{Th sigma}Let $(X_{n})_{n\in Z}$ and $(S_{n})_{n\geq1}$ as defined in
(\ref{def X,S}) and assume $(\xi_{n})$\ has trivial tail sigma field
$\mathcal{T}_{d}$. Then the following are equivalent: \newline(i) $\left(
S_{n}^{2}/n\right)  _{n\geq1}$ is uniformly integrable.\newline(ii) There is
$c\geq0$ such that%
\[
\lim_{n\rightarrow\infty}\frac{E(S_{n}^{2})}{n}=c^{2}\text{ and }\frac{S_{n}%
}{\sqrt{n}}\Rightarrow N(0,c^{2})\text{ as }n\rightarrow\infty.
\]

\end{theorem}

\subsection{Two classes of stationary Markov Chains with trivial
$\mathcal{T}_{d}.$}

\qquad We shall give here two examples of Markov chains with $\mathcal{T}_{d}$ trivial.

\bigskip

\textbf{Absolutely regular Markov chains \label{abs reg}}

For a stationary Markov chain $\bar{\xi}=(\xi_{k})_{k\in Z}$ with values in a
separable Banach space endowed with the Borel sigma algebra $\mathcal{B}$, the
coefficient of absolutely regularity is defined by (see Proposition 3.22 in
Bradley, 2007)
\[
\beta_{n}=\beta_{n}(\mathbf{\xi)=}\beta(\xi_{0},\xi_{n}\mathbf{)=}E\left(
\sup_{A\in\mathcal{B}}|\mathbb{P}(\xi_{n}\in A|\xi_{0})-\mathbb{P}(\xi_{0}\in
A)|\right)  ,
\]
where $\mathcal{B}$ denotes the Borel sigma filed.

Equivalently, (see Corollary 3.30 in Bradley (2007))
\[
\beta_{n}=\beta_{n}(\mathbf{\xi)=}\beta(\xi_{0},\xi_{n}\mathbf{)=}\sup
_{C\in\mathcal{B}^{2}}|P((\xi_{0},\xi_{n})\mathbf{\in}C)-P((\xi_{0},\xi
_{n}^{\ast})\mathbf{\in}C)|,
\]
where $(\xi_{0},\xi_{n}^{\ast}\mathbf{)}$ are independent and identically
distributed. This coefficient was introduced by Volkonskii and Rozanov (1959)
and was attributed there to Kolmogorov.

If $\beta_{n}\rightarrow0,$ the Markov chain is called absolutely regular and
the tail sigma field $\mathcal{T}_{d}$ is trivial (see Section 2.5 in Bradley
(2010)). It follows that both Theorem \ref{ThCLT} and Theorem \ref{Th sigma} hold.

Let us mention that there are numerous examples of stationary absolutely
regular Markov chains. We know that a strictly stationary, countable state
Markov chain is absolutely regular if and only if the chain is irreducible and
aperiodic. Also, any strictly stationary Harris recurrent and aperiodic Markov
chain is absolutely regular. For easy reference we refer to Section 3 in
Bradley (2005) survey paper and to the references mentioned there.

In general, the CLT for this class requires the knowledge of the rates of
convergence to $0$ of the $(\beta_{n})$ coefficients (see for instance Doukhan
et al. (1994) and Peligrad (2020) for a discussion on the CLT under $\beta
_{n}\rightarrow0)$.

\bigskip

\textbf{Interlaced mixing Markov chains}

\qquad Another example where our results apply is the class of interlaced
mixing Markov chains. Let $\mathcal{A},\mathcal{B}$ be two sub $\sigma
$-algebras of $\mathcal{F}$. Define the maximal coefficient of correlation
\[
\rho(\mathcal{A},\mathcal{B})=\sup_{f\in\mathbb{L}_{0}^{2}(\mathcal{A}),\text{
}g\in\mathbb{L}_{0}^{2}(\mathcal{B})}\frac{|E(fg)|}{||f||\cdot||g||}\text{ ,}%
\]
where $\mathbb{L}_{0}^{2}(\mathcal{A})$ ($\mathbb{L}_{0}^{2}(\mathcal{B})$) is
the space of random variables that are $\mathcal{A-}$measurable (respectively
$\mathcal{B-}$measurable$\mathcal{)}$, zero mean and square integrable. For a
sequence of random variables, $(\xi_{k})_{k\in\mathbb{Z}}$, we define
\[
\rho_{n}^{\ast}=\sup\rho(\sigma(\xi_{i},i\in S),\sigma(\xi_{j},j\in
T))\text{,}%
\]
where the supremum is taken over all pairs of disjoint sets, $T$ and $S$ or
$\mathbb{R}$ such that $\min\{|t-s|:t\in T,$ $s\in S\}\geq n.$ We call the
sequence $\rho^{\ast}-$mixing if $\rho_{n}^{\ast}\rightarrow0$ as
$n\rightarrow\infty.$

The $\rho^{\ast}$-mixing condition goes back to Stein (1972) and to Rosenblatt
(1972). It is well-known that $\rho^{\ast}-$mixing implies that the tail sigma
field $\mathcal{T}_{d}$ is trivial (see Section 2.5 in Bradley (2010)). It
follows that both Theorem \ref{ThCLT} and Theorem \ref{Th sigma} hold.
Although these theorems are not new for this class, the results in this paper
provide an unified approach for different classes of Markov chains. For
further reaching results concerning $\rho^{\ast}-$mixing sequences see for
instance Theorem 11.18 in Bradley (2007) and Corollary 9.16 in Merlev\`{e}de,
Peligrad and Utev (2019).

These two classes, absolutely regular and interlaced mixing Markov chains, are
of independent interest. There are known examples (see Example 7.16 in
Bradley, 2007) of $\rho^{\ast}-$mixing sequences which are not absolutely
regular. On the other hand there are known examples of absolutely regular
Markov chains which are not $\rho^{\ast}-$mixing. An example of reversible,
absolutely regular Markov chains which is not $\rho^{\ast}-$mixing was
constructed by Bradley (2015).

\section{Proofs}

\qquad The proofs of both theorems are based on the following result, which is
Theorem 1 in Peligrad (2020), combined with Lemma \ref{Lemma limit} below.

For reader's convenience, let us state first the main result in Peligrad
(2020). It uses the notion of totally ergodic Markov chain. To explain it, let
us consider the operator $P$ induced by the kernel $P(x,A)$ on bounded
measurable functions on $(S,\mathcal{A})$ defined by $Pf(x)=\int
\nolimits_{S}f(y)P(x,dy)$. We call $(\xi_{n})_{n\in\mathbb{Z}}$ totally
ergodic\ if and only if the powers $P^{m}\ $are ergodic with respect to $\pi,$
for all $m\in N$ (i.e. $P^{m}f=f$ for $f$ bounded on $(S,\mathcal{A})$ implies
$f$ is constant $\pi-$a.s.).

\begin{theorem}
\label{Th random center}Let $(X_{n})_{n\in Z}$ and $(S_{n})_{n\geq1}$ as
defined in (\ref{def X,S}), $(\xi_{n})$ is totally ergodic and (\ref{supSn})
is satisfied. Then, the following limit exists
\begin{equation}
\lim_{n\rightarrow\infty}\frac{1}{n}||S_{n}-E(S_{n}|\xi_{0},\xi_{n}%
)||^{2}=c^{2} \label{def teta}%
\end{equation}
and%
\[
\frac{S_{n}-E(S_{n}|\xi_{0},\xi_{n})}{\sqrt{n}}\Rightarrow N(0,c^{2})\text{ as
}n\rightarrow\infty.
\]

\end{theorem}

Next lemma deals with the random centering in Theorem \ref{Th random center}.

\begin{lemma}
\label{Lemma limit}Let $(\xi_{n})_{n\in Z}$ be a stationary sequence not
necessarily Markov with trivial tail sigma field $\mathcal{T}_{d}$. Let
$(X_{n})_{n}$ and $(S_{n})_{n}$ as defined in (\ref{def X,S}). Then
\begin{equation}
E\left\vert E(\frac{S_{n}}{\sqrt{n}}{\LARGE |}\mathcal{F}_{0}\vee
\mathcal{F}^{n})\right\vert \rightarrow0\text{ as }n\rightarrow\infty.
\label{limit L1}%
\end{equation}
If in addition we assume that $\left(  S_{n}^{2}/n\right)  _{n\geq1}$ is
uniformly integrable then%
\begin{equation}
E\left(  E^{2}(\frac{S_{n}}{\sqrt{n}}{\LARGE |}\mathcal{F}_{0}\vee
\mathcal{F}^{n})\right)  \rightarrow0\text{ as }n\rightarrow\infty.
\label{limit L2}%
\end{equation}

\end{lemma}

\textbf{Proof of Lemma \ref{Lemma limit}}.

\bigskip

To prove (\ref{limit L1}) it is clear that it is enough to prove that from any
subsequence of indexes $(n^{\prime})$ convergent to infinity, we can extract
one $(n"),$ also convergent to infinity, and such that (\ref{limit L1}) holds
along $(n")$. Obviously, condition (\ref{supSn}) implies that $(S_{n}/\sqrt
{n})$ is tight. Denote by $\bar{\xi}=(\xi_{n})_{n\in Z}.$ Consider the vector
$W_{n}=(S_{n}/\sqrt{n},\bar{\xi})$ defined in a canonical way on $R\times
S^{Z}$ with values in $R\times S^{Z}$. Note that $W_{n}$ is tight because
$(S_{n}/\sqrt{n})$ is tight and $\bar{\xi}$ does not depend on $n$. Therefore,
from any subsequence $(n^{\prime})$ we can extract one $(n")$ such that
$W_{n"}$ is convergent in distribution, say $W_{n"}\Rightarrow W=(L,\bar{\xi
}^{\prime})$, where $\bar{\xi}^{\prime}$ is distributed as $\bar{\xi}$ and
$S_{n}/\sqrt{n}\Rightarrow L$. Because $R\times S^{Z}$ is separable, by the
Skorohod representation theorem, (see Theorem 6.7 in Billingsley (1999)), we
can expand the probability space to $(\tilde{\Omega},\mathcal{\tilde{F}%
},\tilde{P})$ and construct on this expanded probability space, vectors
$\tilde{W}_{n"}=(\tilde{S}_{n"},\tilde{\xi}^{n"})$ and $\tilde{W}=(\tilde
{L},\tilde{\xi}^{\prime})$ such that for each $n",$ $\tilde{W}_{n"}$ is
distributed as $W_{n"},$  $\tilde{W}\ $is distributed as $W$, and $\tilde
{W}_{n"}\rightarrow\tilde{W}$ a.s.$\,$Note that for each $n"$ we have
$\tilde{\xi}^{n"}=\tilde{\xi}^{\prime}$ a.s. and so $(\tilde{S}_{n"}%
,\tilde{\xi}^{n"})=(\tilde{S}_{n"},\tilde{\xi}^{\prime})$ a.s. Denote
$\mathcal{\tilde{F}}_{n}=\sigma(\tilde{\xi}_{k}^{\prime},k\leq n)$ and
$\mathcal{\tilde{F}}^{n}=\sigma(\tilde{\xi}_{k}^{\prime},k\geq n),$ completed
with the sets of measure $0$ and $\mathcal{\tilde{T}}_{d}=\cap_{n\geq
1}(\mathcal{\tilde{F}}_{-n}\vee\mathcal{\tilde{F}}^{n})$. Note that the
Skorohod representation (see page 71 in Billingsley (1999)) starts with the
construction of $(\tilde{L},\tilde{\xi}^{\prime})$ in a canonical way on
$R\times S^{Z}$, such that the marginals are distributed as $L$ and $\bar{\xi
}$. But because $\bar{\xi}$ was also constructed in a canonical way, we have
that the tail sigma field $\mathcal{\tilde{T}}_{d}$ of $\tilde{\xi}^{\prime}$
is also trivial.

To simplify the notation let us re-denote the index $n"$ by $n$. Clearly,%

\[
\frac{\tilde{S}_{n}}{\sqrt{n}}\rightarrow\tilde{L}\text{ a.s. as }%
n\rightarrow\infty.
\]
Now (\ref{supSn}) implies that $(\tilde{S}_{n}/\sqrt{n})$ is uniformly
integrable, so we also have%
\begin{equation}
\tilde{E}\left\vert \frac{\tilde{S}_{n}}{\sqrt{n}}-\tilde{L}\right\vert
\rightarrow0\text{ as }n\rightarrow\infty, \label{weakL1}%
\end{equation}
and because $E(X_{1})=0$, by the convergence of moments in the weak laws (see
Theorem 3.5 in Billingsley (1999)), we have that
\begin{equation}
\tilde{E}(\tilde{L})=0. \label{E(L)0}%
\end{equation}
By the Fatou lemma, we also have that $\tilde{E}(\tilde{L}^{2})<\infty.$

By stationarity and the triangle inequality, note that for every $m\in N,$
$m\leq n,$%
\begin{gather}
E\left\vert E\left(  \frac{S_{n}}{\sqrt{n}}|\mathcal{F}_{0}\vee\mathcal{F}%
^{n}\right)  \right\vert \leq E\left\vert E\left(  \frac{S_{n}-S_{m}}{\sqrt
{n}}|\mathcal{F}_{0}\vee\mathcal{F}^{n}\right)  \right\vert + \label{equality}%
\\
E\left\vert E\left(  \frac{S_{n}-S_{m}}{\sqrt{n}}|\mathcal{F}_{0}%
\vee\mathcal{F}^{n}\right)  -E\left(  \frac{S_{n}}{\sqrt{n}}|\mathcal{F}%
_{0}\vee\mathcal{F}^{n}\right)  \right\vert \nonumber\\
\leq E\left\vert E\left(  \frac{S_{n-m}}{\sqrt{n}}{\LARGE |}\mathcal{F}%
_{-m}\vee\mathcal{F}^{n-m}\right)  \right\vert +\frac{E|S_{m}|}{\sqrt{n}%
}.\nonumber
\end{gather}
Because \thinspace$W_{n},$ and $\tilde{W}_{n"}$ have the same distribution,
\begin{equation}
E\left\vert E\left(  \frac{S_{n-m}}{\sqrt{n}}{\LARGE |}\mathcal{F}_{-m}%
\vee\mathcal{F}^{n-m}\right)  \right\vert =\tilde{E}\left\vert \tilde
{E}\left(  \frac{\tilde{S}_{n-m}}{\sqrt{n}}{\LARGE |}\mathcal{\tilde{F}}%
_{-m}\vee\mathcal{\tilde{F}}^{n-m}\right)  \right\vert . \label{identity}%
\end{equation}
Now we use the following inequality: $\ $%
\begin{gather}
\tilde{E}\left\vert \tilde{E}\left(  \frac{\tilde{S}_{n-m}}{\sqrt{n}%
}{\LARGE |}\mathcal{\tilde{F}}_{-m}\vee\mathcal{\tilde{F}}^{n-m}\right)
\right\vert \label{decompose}\\
\leq\tilde{E}\left\vert \tilde{E}\left(  \frac{\tilde{S}_{n-m}}{\sqrt{n}%
}-\tilde{L}{\LARGE |}\mathcal{\tilde{F}}_{-m}\vee\mathcal{\tilde{F}}%
^{n-m}\right)  \right\vert +\tilde{E}{\LARGE |}\tilde{E}{\LARGE (}\tilde
{L}{\LARGE |}\mathcal{\tilde{F}}_{-m}\vee\mathcal{\tilde{F}}^{n-m}%
{\LARGE )|}.\nonumber
\end{gather}
We treat now the first term in left hand side of (\ref{decompose}). By the
properties of the conditional expectation
\begin{equation}
\tilde{E}\left\vert \tilde{E}\left(  \frac{\tilde{S}_{n-m}}{\sqrt{n}}%
-\tilde{L}{\LARGE |}\mathcal{\tilde{F}}_{-m}\vee\mathcal{\tilde{F}}%
^{n-m}\right)  \right\vert \leq\tilde{E}\left\vert \frac{\tilde{S}_{n-m}%
}{\sqrt{n}}-\tilde{L}\right\vert . \label{contract}%
\end{equation}
Overall, starting from (\ref{equality}) combined to (\ref{identity}),
(\ref{decompose}) and (\ref{contract}), we obtain for $n,m\in N,m\leq n,$
\begin{gather*}
E\left\vert E\left(  \frac{S_{n}}{\sqrt{n}}|\mathcal{F}_{0}\vee\mathcal{F}%
^{n}\right)  \right\vert \leq\\
\tilde{E}\left\vert \frac{\tilde{S}_{n-m}}{\sqrt{n}}-\tilde{L}\right\vert
+\tilde{E}\left\vert \tilde{E}{\LARGE (}\tilde{L}{\LARGE |}\mathcal{\tilde{F}%
}_{-m}\vee\mathcal{\tilde{F}}^{n-m}{\LARGE )}\right\vert +\frac{E|S_{m}%
|}{\sqrt{n}}.
\end{gather*}
Therefore, for $m\in N$ fixed, by letting $n\rightarrow\infty,$ and by taking
into account (\ref{weakL1})\ and the fact that $\mathcal{\tilde{F}}_{-m}%
\vee\mathcal{\tilde{F}}^{n-m}$ is decreasing in $n$
\begin{equation}
\lim\sup_{n}E\left\vert E\left(  \frac{S_{n}}{\sqrt{n}}|\mathcal{F}_{0}%
\vee\mathcal{F}^{n}\right)  \right\vert \leq\tilde{E}\left\vert \tilde
{E}{\LARGE (}\tilde{L}{\LARGE |}\mathcal{\cap}_{n\geq1}\left(  \mathcal{\tilde
{F}}_{-m}\vee\mathcal{\tilde{F}}^{n}\right)  )\right\vert . \label{limit n}%
\end{equation}
Now, by letting \thinspace$m\rightarrow\infty$, and using the fact that
$\mathcal{\cap}_{n\geq1}\left(  \mathcal{\tilde{F}}_{-m}\vee\mathcal{\tilde
{F}}^{n}\right)  $ is decreasing in $m,$ we obtain by (\ref{E(L)0}) that
\begin{align}
\lim\sup_{n}E\left\vert E\left(  \frac{S_{n}}{\sqrt{n}}|\mathcal{F}_{0}%
\vee\mathcal{F}^{n}\right)  \right\vert  &  \leq\tilde{E}\left\vert \tilde
{E}(\tilde{L}|\mathcal{\cap}_{m\geq1}\mathcal{\cap}_{n\geq1}\left(
\mathcal{\tilde{F}}_{-m}\vee\mathcal{\tilde{F}}^{n}\right)  )\right\vert
\label{limit n and m}\\
&  \leq\tilde{E}\left\vert \tilde{E}(\tilde{L}|\mathcal{\tilde{T}}%
_{d})\right\vert =\left\vert \tilde{E}(\tilde{L})\right\vert =0.\nonumber
\end{align}

$\square$

\bigskip

\textbf{Proof of Theorem \ref{ThCLT}}

\bigskip

First of all we mention that, by Proposition 2.12 in the Vol. 1 of Bradley
(2007), the sequence $(\xi_{k})_{k\in Z}$ is totally ergodic. Since we assumed
(\ref{supSn}), by Theorem \ref{Th random center}
\[
\frac{S_{n}-E(S_{n}|\xi_{0},\xi_{n})}{\sqrt{n}}\Rightarrow N(0,c^{2})\text{ as
}n\rightarrow\infty.
\]
By Lemma \ref{Lemma limit}, and by the Markov property
\[
E\left\vert E(\frac{S_{n}}{\sqrt{n}}{\LARGE |}\xi_{0},\xi_{n})\right\vert
=\left\vert E(\frac{S_{n}}{\sqrt{n}}|\mathcal{F}_{0}\vee\mathcal{F}%
^{n})\right\vert \rightarrow0\text{ as }n\rightarrow\infty.
\]
Therefore
\[
\frac{S_{n}}{\sqrt{n}}\Rightarrow N(0,c^{2})\text{ as }n\rightarrow\infty.
\]
Now, because $E(S_{n}^{2}/n)<\infty$ it follows that (see Theorem 3.5 in
Billingsley (1999))
\[
\lim_{n\rightarrow\infty}E\frac{|S_{n}|}{\sqrt{n}}=E|N(0,c^{2})|=\frac
{c}{\sqrt{2\pi}}.
\]
So, $c$ can be identified as
\[
c=\sqrt{2\pi}\lim_{n\rightarrow\infty}E\frac{|S_{n}|}{\sqrt{n}}.
\]
$\square$

\textbf{Proof of Theorem \ref{Th sigma}.}

Assume (i). Since $\left(  S_{n}^{2}/n\right)  _{n\geq1}$ is uniformly
integrable it follows that Theorem \ref{ThCLT} holds. In addition, (see
Theorem 3.5 in Billingsley (1999))
\[
\lim_{n\rightarrow\infty}\frac{E(S_{n}^{2})}{n}=c^{2}.
\]
Because by Theorem \ref{Th random center}%

\[
\lim_{n\rightarrow\infty}\frac{1}{n}||S_{n}-E(S_{n}|\xi_{0},\xi_{n}%
)||^{2}=c^{2}=\lim_{n\rightarrow\infty}\frac{1}{n}||S_{n}||^{2}-\lim
_{n\rightarrow\infty}\frac{1}{n}||E(S_{n}|\xi_{0},\xi_{n})||^{2},
\]
we also have
\[
\lim_{n\rightarrow\infty}\frac{1}{n}||E(S_{n}|\xi_{0},\xi_{n})||^{2}=0.
\]
and so (i) implies (ii). On the other hand, if (ii) holds, by the convergence
of moments in the CLT (Theorem 3.6 in Billingsley (1999)), we have that
$\left(  S_{n}^{2}/n\right)  _{n\geq1}$ is uniformly integrable. $\ \square$

\textbf{Acknowledgement.} The author would like to thank to Yizao Wang for
useful discussions.

\end{document}